\documentclass[11pt]{amsart}

\usepackage{amsfonts,amssymb,amscd,amsmath,latexsym,amsbsy,amsthm}

\usepackage{amssymb}
\usepackage{amsfonts}
\usepackage{latexsym}

\theoremstyle{plain}
\newtheorem{theorem}{Theorem}
\newtheorem{thm}[theorem]{Theorem}

\newtheorem{lem}[theorem]{Lemma}

\newtheorem{prop}[theorem]{Proposition}

\newtheorem{cor}[theorem]{Corollary}

\theoremstyle{definition}

\newtheorem*{defn}{Definition}

\theoremstyle{remark}
\newtheorem*{remark}{Remark}
\newtheorem*{rems}{Remarks}
\newtheorem*{rem}{Remark}

\numberwithin{theorem}{section}

\newcommand{\Fin}{\text{\bf Fin}}
\newcommand{\Bij}{\text{\bf Bij}}
\newcommand{\GrAlg}{{\QQ\text{-\bf GrAlg}}}
\newcommand{\GrVect}{{\QQ\text{-\bf GrVect}}}

\DeclareMathOperator{\Exp}{Exp}
\DeclareMathOperator{\Log}{Log}
\DeclareMathOperator{\Sinh}{Sinh}
\DeclareMathOperator{\Arcsinh}{Arcsinh}
\DeclareMathOperator{\Cosh}{Cosh}
\DeclareMathOperator{\sgn}{sgn}
\DeclareMathOperator{\Sym}{Sym}
\DeclareMathOperator{\Tr}{Tr}
\DeclareMathOperator{\arcsinh}{arcsinh}
\DeclareMathOperator{\odd}{odd}

\newcommand{\2}{_{(2)}}

\newcommand{\ben}{\begin{enumerate}}
\newcommand{\een}{\end{enumerate}}

\newcommand{\QQ}{{\mathbb{Q}}}
\newcommand{\CC}{{\mathbb{C}}}
\newcommand{\RR}{{\mathbb{R}}}

\newcommand{\ZZ}{{\mathbb{Z}}}

\newcommand{\solu}[1]{\begin{sol}{\bf (\ref{#1})}}

\begin{document}

\title{The action of $S_n$ on the cohomology of $\overline{M_{0,n}}(\RR)$}

\author{Eric Rains}

\begin{abstract}
In recent work by Etingof, Henriques, Kamnitzer, and the author, a
presentation and explicit basis was given for the rational cohomology of
the real locus $\overline{M_{0,n}}(\RR)$ of the moduli space of stable
genus 0 curves with $n$ marked points.  We determine the graded character
of the action of $S_n$ on this space (induced by permutations of the marked
points), both in the form of a plethystic formula for the cycle index, and
as an explicit product formula for the value of the character on a given
cycle type.
\end{abstract}

\maketitle

\section{Introduction}

For any integer $n\ge 3$, let $\overline{M_{0,n}}$ be the moduli space of
stable curves of genus 0 with $n$ marked points; by convention, for $n=1$,
$n=2$, this is just a single point, but we never allow $n=0$.  Since a
stable curve of genus 0 has trivial automorphism group, this is in fact a
smooth projective scheme over $\ZZ$ (and a {\em fine} moduli space), and
thus its real locus $M_n:=\overline{M_{0,n}}(\RR)$ is a smooth compact
manifold.  The symmetric group acts on $M_n$ by permuting the marked
points, and thus acts on the cohomology.  The main result of the present
work is an explicit product formula for the (graded) character of this
action.

\begin{thm}\label{thm:main}
Let $\pi\in S_n$ be a permutation with $n_1+1$ fixed points and $n_m$
$m$-cycles for $m>1$, and define
\[
o_m = \sum_{1\le k:2^k|m} (2^{-k}m) n_{2^{-k}m}.
\]
Then
\begin{align}
\sum_k &(-t)^k \Tr(\pi|H^k(M_n,\QQ))\notag\\
&=
\prod_{1\le l}
(\gamma_l(t)+o_l t^{l/2})
\prod_{0\le i\le n_l-2} (\gamma_l(t)+(o_l+l(n_l-2-2i))t^{l/2}),
\notag
\end{align}
where the polynomials $\gamma_l(t)$ satisfy
\[
\sum_{\text{odd }k|l} t^{-l/2k}\gamma_{l/k}(t) = t^{-l/2}.
\]
\end{thm}

\begin{rems}
1. Note that we are using the standard convention for products with negatively
many terms; thus for $n_l\le 0$,
\begin{align}
\prod_{0\le i\le n_l-2} &(\gamma_l(t)+(o_l+l(n_l-2-2i))t^{l/2})\notag\\
&{}:=
\prod_{n_l-1\le i\le -1} (\gamma_l(t)+(o_l+l(n_l-2-2i))t^{l/2})^{-1}.
\notag
\end{align}
In particular, the infinite product is indeed well-defined, since if
$n_l=0$, the corresponding factor is 1.  For $n_1=-1$, the corresponding
factor is
\[
\prod_{-2\le i\le -1} (1+(-3-2i)t^{1/2})^{-1}
=
1/(1-t).
\]
Similarly, the presence of $t^{l/2}$ for $l$ odd is not an issue, since
then the corresponding factor is invariant under $t^{l/2}\to -t^{l/2}$
(simply reverse the order of multiplication in the product over $i$, and
note that $o_l=0$).  Finally, $\gamma_l(t)$ is indeed a polynomial, since
by M\"obius inversion,
\[
t^{-l/2}\gamma_l(t) = \sum_{\text{odd }k|l} \mu(k)t^{-l/2k},
\]
and thus
\[
\gamma_l(t) = \sum_{\text{odd }k|l} \mu(k) t^{l(1-1/k)/2}.
\]

2. Also note the factor $(-1)^k$ above; in particular, the Euler character
of $M_n$ is given by setting $t=1$ above (or taking a limit, if $n_1=-1$).
In this context, it is worth noting that $\gamma_l(1)=0$ unless $l$ is a
power of 2, and $\gamma_{2^k}(t)=1$.

3.
On the identity element, we obtain
\[
\prod_{0\le i\le n-3} (1 + (n-3-2i) t^{1/2})
=
\prod_{0\le i\le \lfloor(n-3)/2\rfloor} (1-(n-3-2i)^2 t),
\]
agreeing with the formula of \cite{EHKR} for the Poincar\'e series of
$M_n$.

4. We finally note that the above formula is remarkably similar to the
following formula of Lehrer \cite{Leh,KL}, valid for $n\ge 3$:
\[
\sum_k (-t)^k \Tr(\pi|H^k(M_{0,n}(\CC),\QQ))
=
(1-t)^{-1}
\prod_{1\le l}
\prod_{0\le i<n_l} (\eta_l(t)-lit^l),
\]
where the polynomials $\eta_l(t)$ satisfy
\[
\sum_{k|l} t^{-l/k}\eta_l(t) = t^{-l}.
\]
Of course, the close analogy between the cohomology of these spaces was
already noted in \cite{EHKR}.
\end{rems}

As one might imagine from the form of the above result, it is much more
natural to consider the action of $S_{n-1}$ on $M_n$, rather than the full
action of $S_n$.  Indeed, the results of \cite{EHKR} on the structure of
$H^*(M_n,\QQ)$ (summarized in Section 2) give a particularly nice
description of this restriction in terms of the homology (not {\em
  co}homology, as one would normally expect) of a certain poset; the
corresponding character was studied in \cite{CHR}.  In Section 3, by
combining these results, we obtain an expression (Theorem \ref{thm:cycle})
for the ``cycle index'' of the restriction, i.e., a generating function for
the character.  In Section 4, we derive a number of differential equations
satisfied by the cycle index; the corresponding recurrences for the
character prove the theorem for the restriction (i.e., when $\pi$ has a
fixed point).  Finally, in Section 5, we show that $H^*(M_n,\QQ)$ satisfies
a particularly strong form of functoriality which in particular enables us
to derive the full $S_n$ character from the $S_{n-1}$ character alone,
proving the main theorem.  (We also give an expression for the
corresponding cycle index (Theorem \ref{thm:ext_cycle}).)  Finally, in
Corollary \ref{cor:euler}, we give a formula for the Euler character of
$M_n$, in particular determining the precise permutations for which the
Euler character is nonzero.

\noindent{\bf Notation}

As we are dealing with cohomology, it will be convenient to use ``super''
conventions.  That is, if $V_1$,\dots $V_n$ is a sequence of graded vector
spaces (with it being understood here and in the sequel that the
coefficient field is $\QQ$ and all nontrivial homogeneous components have
finite dimension and nonnegative degree), we identify the two tensor
products
\[
V_1\otimes V_2\otimes\cdots\otimes V_n
\]
and
\[
V_{\pi(1)}\otimes V_{\pi(2)}\otimes\cdots\otimes V_{\pi(n)}
\]
for any permutation $\pi$ via the isomorphism
\[
v_1\otimes v_2\otimes\cdots\otimes v_n
\to
\prod_{i<j,\pi(i)>\pi(j)} (-1)^{\deg(v_i)\deg(v_j)}
v_{\pi(1)}\otimes v_{\pi(2)}\otimes\cdots\otimes v_{\pi(n)}
\]
for any sequence of homogeneous elements $v_i\in V_i$.  Similarly, if $A$
is a graded algebra, we say that it is supercommutative if
\[
xy = (-1)^{\deg(x)\deg(y)} yx.
\]
In particular, the free supercommutative algebra generated by elements of
degree 1 is simply the exterior algebra.

\noindent{\bf Acknowledgements}

The author would like to thank P. Etingof for useful discussions (and for
asking the question in the first place); also A. Henderson for helpful
comments and references regarding the complex case.  This work was
supported in part by NSF Grant No. DMS-0401387.

\section{The cohomology of $M_n$}

\begin{thm}\cite{EHKR}
For $n\ge 1$, the algebra $\Lambda_n:=H^*(M_n,\QQ)$ is the supercommutative
quadratic algebra generated over $\QQ$ by elements $\omega_{ijkl}, 1\le
i,j,k,l\le n$, antisymmetric in $ijkl$, with defining relations
\[
\omega_{ijkl}+\omega_{jklm}+\omega_{klmi}+\omega_{lmij}+\omega_{mijk}=0
\]
and
\[
\omega_{ijkl}\omega_{ijkm}
\]
for any distinct $i,j,k,l,m$.  Moreover, the action of $S_n$ on
$H^*(M_n,\QQ)$ is given in terms of these generators by
\[
\pi^*(\omega_{ijkl}) =
\omega_{\pi(i)\pi(j)\pi(k)\pi(l)}.
\]
\end{thm}

This extends naturally to a functor $\Lambda:\Bij^+\to\GrAlg$, where
$\Bij^+$ is the category of nonempty finite sets and bijections, and
$\GrAlg$ is the category of graded $\QQ$-algebras.  As we mentioned in
the introduction, we will need to consider also a restriction of this to
the category $\Bij$ of all finite sets and bijections.

\begin{prop}\cite{EHKR}
For any ordered finite set $S$, let $\Lambda'(S)$ denote the
supercommutative algebra generated by antisymmetric elements $\nu_{ijk}$
for distinct $i,j,k\in S$ subject to the relations
\[
\nu_{ijk}\nu_{ijl} = 0
\]
and
\[
\nu_{ijk}\nu_{klm}+\nu_{jkl}\nu_{lmi}+\nu_{klm}\nu_{mij}+\nu_{lmi}\nu_{ijk}+\nu_{mij}\nu_{jkl}
=
0;
\]
extend this to a functor $\Bij\to \GrAlg$ by
\[
\Lambda'(\pi)(\nu_{ijk})=\nu_{\pi(i)\pi(j)\pi(k)}.
\]
Then for each $n\ge 0$, there is an isomorphism
$\Lambda'(\{1,2,\dots,n\})\cong \Lambda_{n+1}$ defined on generators by
\[
\nu_{ijk}\mapsto \omega_{ijkn}.
\]
\end{prop}

A monomial in the generators $\nu_{ijk}$ determines an equivalence relation
on $S$ (taking $i\cong j\cong k$ if $\nu_{ijk}$ appears in the monomial);
equivalently, each monomial determines a partition of $S$ into (unordered)
disjoint subsets.  If $\rho$ is such a partition (a fact denoted by the
relation $\rho\vdash S$), let $\Lambda'[\rho]$ denote the span in
$\Lambda'(S)$ of all monomials corresponding to $\rho$; note that
$\Lambda'[\rho]$ is unchanged (up to canonical isomorphism) if we remove a
singleton class from $\rho$ and $S$.  In particular, we may let
$\Lambda'[T]$ denote the case in which $\rho$ has a single nontrivial
equivalence class, equal to $T$; the result is independent of $S$ up to
canonical isomorphism.

\begin{thm}\cite{EHKR}
The spaces $\Lambda'[\rho]$ for different $\rho$ are linearly independent,
and thus
\[
\Lambda'(S) = \bigoplus_{\rho\vdash S} \Lambda'[\rho].
\]
If $\rho$ has classes
$\rho_1$, $\rho_2$,\dots, $\rho_k$, then multiplication in $\Lambda'(S)$
induces a natural isomorphism
\[
\Lambda'[\rho_1,\rho_2,\dots,\rho_k] \cong \Lambda'[\rho_1]\otimes
\Lambda'[\rho_2]\otimes\cdots\otimes\Lambda'[\rho_k];
\]
this remains valid even if some singleton classes of $\rho$ are omitted.
\end{thm}

Finally, the indecomposable spaces $\Lambda'[T]$ can be expressed in terms
of certain poset homology groups.

\begin{thm}\cite{EHKR}
If $|T|$ is even, then $\Lambda'[T]=0$; otherwise, if $|T|=2n+1$,
\[
\Lambda'[T] \cong \tilde{H}_n(\Pi^{\odd}_{T},\QQ)\otimes \sgn,
\]
where $\Pi^{\odd}_{T}$ is the poset of partitions of $T$ with all parts
odd, $\tilde{H}_n$ is the top (shifted) reduced homology of this poset, and
$\sgn$ is the sign representation of $\Sym(T)$.
\end{thm}

\begin{rem}
Note that in $\tilde{H}_n$, the degree has been shifted by 1 from the
standard definition of poset homology, in order to obtain the correct
degree in $\Lambda'[T]$.  In any event, $\tilde{H}_n$ is the only
nontrivial homology group of $\Pi^{\odd}_T$ (which is Cohen-Macaulay
\cite{Bj,CHR}), so there is no risk of confusion.
\end{rem}

\section{Cycle indices}

Let $\GrVect$ denote the category of graded vector spaces $W$ and
degree 0 linear transformations.  Given an endomorphism $\phi:W\to W$ in
$\GrVect$, the graded trace of $\phi$ is defined to be the power
series $\Tr(\phi)\in \QQ[[t]]$ defined by
\[
\Tr(\phi)(t):=\sum_{k\ge 0} t^k (-1)^k \Tr_{W_k}(\phi);
\]
the sign factor reflects our interpretation of $W$ as a graded superspace.

Now, let $V$ be a representation of $\Bij$ in $\GrVect$ (a
``graded representation of $\Bij$'').

\begin{defn}
The {\em cycle index} of $V$ is the power series $Z_V\in
\QQ[[t,p_1,p_2,\dots]]$ given by
\[
\sum_{n\ge 0}
\frac{1}{n!}
\sum_{\pi\in S_n} \Tr(V(\pi))(t) \prod_i p_i^{n_i(\pi)}
\]
where for a permutation $\pi$, $n_i(\pi)$ is the number of $i$-cycles of
$\pi$.
\end{defn}

\begin{rem}
We may similarly associate a cycle index to an arbitrary virtual (graded)
character of $\Bij$ (i.e., a sequence $\chi_n$ such that $\chi_n$ is
a virtual character of $S_n$).
\end{rem}

There are two natural gradings on the above algebra of power series
($t$-degree and $p$-degree), defined on generators by
\[
\deg_t(t)=1,\deg_t(p_i)=0,
\deg_p(t)=0,\deg_p(p_i)=i;
\]
a cycle index $Z_V$ is homogeneous of $t$-degree $d$ if $V$ is homogeneous
of degree $d$, and homogeneous of $p$-degree $d$ if $V(S)=0$ for $|S|\ne
d$.

The sum and product of cycle indices is itself a cycle index, as is
\[
F^\sim:=F(t,p_1,-p_2,p_3,-p_4,\dots).
\]

\begin{prop}
Let $V$ and $W$ be two graded representations of $\Bij$.
Then
\[
Z_V + Z_W = Z_{V\oplus W},
\qquad Z_V Z_W = Z_{V\cdot W},
\qquad Z_V^\sim = Z_{V\otimes \sgn}
\]
where
\begin{align}
(V\oplus W)(S) &= V(S)\oplus W(S),\notag\\
(V\cdot W)(S) &= \bigoplus_{T\subset S} V(T)\otimes W(S\setminus T),\notag
\end{align}
extended to functors in the natural way.
\end{prop}

There is a further operation known as plethysm (or composition), which on
two series $F$ and $G$ with $G(t,0,0,\dots)=0$ is defined as
\[
F[G] := F(t,G(t,p_1,p_2,\dots),G(t^2,p_2,p_4,\dots),\dots);
\]
this is easily verified to be an associative (but not commutative or
distributive) operation.  We will also need the obvious extension of this
to series involving fractional powers of $t$.

\begin{prop}
For any graded representations $V$ and $W$ of $\Bij$ such that
$W(\emptyset)=0$, we have
\[
Z_V[Z_W] = Z_{V[W]},
\]
where $V[W]$ is the graded representation with
\[
V[W](S) := \bigoplus_{\rho\vdash S} V(\rho)\otimes \bigotimes_i W(\rho_i),
\]
extended in the natural way to a functor.
\end{prop}

\begin{rem}
If $W$ is supported on sets of a given cardinality, this is essentially
classical; for the general case, see for instance \cite[Thm.
6.5]{Hen}.
\end{rem}

In general, plethysm does not interact well with tensoring with the sign
character; there is, however, one important special case.

\begin{prop}
If every term of the series $G$ has odd $p$-degree, then
\[
F^\sim[G^\sim] = F[G]^\sim.
\]
\end{prop}

There are three particularly important cycle indices.  For the trivial
representation, we have
\[
\Exp := Z_{\text{triv}} = \exp(\sum_{i\ge 1} p_i/i).
\]
In particular, $\Exp[Z_V]$ is the cycle index of the functor
\[
S\mapsto \bigoplus_{\rho\vdash S} \bigotimes_i V(\rho_i).
\]
We will also need analogues of the hyperbolic sine and cosine:
\begin{align}
\Cosh &:= \frac{\exp(\sum_{i\ge 1} p_i/i)+\exp(\sum_{i\ge 1}(-1)^ip_i/i)}{2}\notag\\
\Sinh &:= \frac{\exp(\sum_{i\ge 1} p_i/i)-\exp(\sum_{i\ge 1}(-1)^ip_i/i)}{2}.\notag
\end{align}
The corresponding representations are obtained from the trivial
representation by removing the spaces associated to sets with odd or even
cardinality, respectively.

We can now state Calderbank, Hanlon, and Robinson's result on the homology
of $\Pi^{\odd}_n$.

\begin{thm}\label{thm:CHR}\cite{CHR}
The cycle index of the functor $\tilde{H}_*(\Pi^{\odd}_{T},\QQ)$ is
\[
(1-\Cosh[\Arcsinh[t^{1/2}p_1]]) + (t^{-1/2} \Arcsinh[t^{1/2}p_1]) ,
\]
where $\Arcsinh$ is the unique symmetric function such that
\[
\Sinh[\Arcsinh]=\Arcsinh[\Sinh]=p_1.
\]
\end{thm}

Note that the first term gives the cycle index for $|T|$ even, while the
second term gives the cycle index for $|T|$ odd.  Also, since $\Sinh$ is
concentrated in odd $p$-degree, the same is true of $\Arcsinh$, and thus
\[
\Sinh^\sim[\Arcsinh^\sim]=\Arcsinh^\sim[\Sinh^\sim]=p_1^\sim=p_1.
\]

This then gives us our first result on the action of $S_n$ on
$H^*(M_n,\QQ)$.

\begin{thm}\label{thm:cycle}
The cycle index of the functor $\Lambda'$ is
\[
\Exp[t^{-1/2}\Arcsinh^\sim[t^{1/2}p_1]]
\]
where $\Arcsinh^\sim$ is the unique symmetric function such that
\[
\Sinh^\sim[\Arcsinh^\sim] = \Arcsinh^\sim[\Sinh^\sim]=p_1.
\]
\end{thm}

\begin{proof}
Since
\[
\Lambda'(S)
\cong
\bigoplus_{\rho\vdash S}
\Lambda'[\rho]
\cong
\bigoplus_{\rho\vdash S}
\bigotimes_i \Lambda'[\rho_i]
\cong
\bigoplus_{\substack{\rho\vdash S\\ \text{all parts odd}}}
\bigotimes_i (\tilde{H}_n(\Pi^{\odd}_S,\QQ)\otimes \sgn),
\]
it follows that
$Z_{\Lambda'} = \Exp[Z_V]$,
where for $|S|$ odd,
\[
V(S)=\tilde{H}_n(\Pi^{\odd}_S,\QQ)\otimes \sgn=\tilde{H}_*(\Pi^{\odd}_S,\QQ)\otimes \sgn.
\]
But tensoring with $\sgn$ simply applies the homomorphism $p_i\to
(-1)^{i-1}p_i$ to the cycle index; the result follows.
\end{proof}

It seems appropriate to mention in passing the corresponding formula for
the cohomology of the {\em complex} moduli space.

\begin{thm}\label{thm:cyc_ind_C} \cite{GK,Hen}
The cycle index of the functor $S\mapsto
H^*(\overline{M_{0,|S|+1}}(\CC),\QQ)$ is given by $\Exp[C]$ where $C$ is
the plethystic inverse of
\[
\frac{t^{-2}(\Exp[t^2 p_1]-1)-t^2(\Exp-1)}{1-t^2}.
\]
\end{thm}

\begin{proof}
The argument of Theorem 4.5 of \cite{EHKR} for computing the Poincar\'e
series from (essentially) the $S_{n-1}$-invariant basis of \cite{Yuz}
extends immediately to the level of cycle indices.  We thus find that the
cycle index is of the form $\Exp[C]$ where $C$ satisfies
\[
C = p_1 + \sum_{m\ge 3} \sum_{1\le l\le m-2} t^{2l} h_m[C];
\]
here $h_m$ is the cycle index of the trivial representation of $S_m$, or
equivalently the $p$-degree $m$ component of $\Exp$.  Moving the sum to the
left-hand side, we find that this indeed specified $C$ as a plethystic
inverse; simplifying the geometric sum gives the desired result.
\end{proof}

\begin{rems}
1. It does not appear to be feasible to obtain an explicit formula for the
graded character (unlike the real case, as we will shortly see); indeed, it
appears that no formula is known for the Poincar\'e series, let alone any
other values of the character.

2. In the references, this is expressed as 1 plus the plethystic inverse of
\[
\Exp[t^2\Log[1+p_1]]/(t^4-t^2) -1/(t^4-t^2) - p_1/(t^2-1),
\]
where $\Log[1+p_1]$ denotes the plethystic inverse of $\Exp-1$.
This in turn is essentially the cycle index of the cohomology of the {\em
  un}-compactified moduli space $M_{0,|S|+1}(\CC)$, suggesting that there
should be a cohomological interpretation for the plethystic inverse
\[
t^{-1/2}\Sinh^\sim[t^{1/2}\Log[1+p_1]]
\]
of the cycle index for $H^*(M_{|S|+1})$.
\end{rems}

\section{The explicit graded character}

It turns out that by using some ideas from \cite{CHR}, we can actually
obtain an explicit formula for the graded character, rather than a mere
generating function.  It will be convenient to introduct another symmetric
function
\[
X := \sum_{k\ \text{odd}} \Arcsinh^\sim[p_k]/k.
\]
It follows from the definition of $\Arcsinh$ that $X$ is the plethystic
inverse of the series
\[
\frac{\exp(p_1-\sum_{k>0}
  p_{2^k}/2^k)-\exp(-p_1-\sum_{k>0}p_{2^k}/2^k)}{2},
\]
and thus is a function of $p_1$, $p_2$, $p_4$,\dots alone (since that
subalgebra is closed under plethysm).  Similarly,
\begin{align}
C:={}&\Cosh^\sim[\Arcsinh^\sim]\notag\\
{}={}&\frac{\exp(X-\sum_{k>0} X[p_{2^k}]/2^k)
+
\exp(-X-\sum_{k>0} X[p_{2^k}]/2^k)}{2}
\notag
\end{align}
also depends only on the variables $p_{2^k}$.

\begin{lem}
The function $X$ satisfies the differential equation
\[
2^l C[p_{2^l}] \frac{\partial}{\partial p_{2^l}} X
=
\delta_{l0}
+
\sum_{0\le k<l} 2^k p_{2^k} \frac{\partial}{\partial p_{2^k}} X.
\]
\end{lem}

\begin{proof}
If we differentiate the plethystic equation
\[
\frac{\exp(X-\sum_{k>0}
  X[p_{2^k}]/2^k)-\exp(-X-\sum_{k>0}X[p_{2^k}]/2^k)}{2},
\]
we find
\[
2^l C \frac{\partial}{\partial p_{2^l}} X
=
\delta_{l0}
+
p_1
\sum_{0\le k<l} (2^k\frac{\partial}{\partial p_{2^k}} X)[p_{2^{l-k}}],
\]
so in particular the claim holds for $l=0$.  For $l>0$, if we multiply both
sides by $C[p_{2^l}]/C$, we find by induction that
{\allowdisplaybreaks
\begin{align}
2^l C[p_{2^l}] \frac{\partial}{\partial p_{2^l}} X
&=
\frac{p_1}{C}
\sum_{0\le k<l} (2^kC[p_{2^k}]\frac{\partial}{\partial p_{2^k}} X)[p_{2^{l-k}}]
\notag\\
&=
\frac{p_1}{C}
+
\sum_{0\le j<k<l} 
(2^j p_{2^j} \frac{\partial}{\partial p_{2^j}} X)[p_{2^{l-k}}]
\notag\\
&=
\frac{p_1}{C}
+
\sum_{0\le j<k<l} 
(2^j p_{2^j} \frac{\partial}{\partial p_{2^j}} X)[p_{2^{k-j}}]
\notag\\
&=
\frac{p_1}{C}
+
\sum_{0<k<l}
p_{2^k}
2^k \frac{\partial}{\partial p_{2^k}} X.\notag
\end{align}
}
\end{proof}

\begin{lem}
Let $c_1$, $c_2$, $c_3$,\dots be indeterminates, and define a symmetric
function
\[
G:=\exp(\sum_{k\ge 0} c_k X[p_k]/k).
\]
Let $G_l$ be the result of setting $p_k=0$ for $k>l$.  Then $G_l$ satisfies
the differential equations
\[
\Bigl(
c_l
+
\sum_{\substack{1\le k\\ 2^k|l}}
2^{-k}l p_{2^{-k}l} \frac{\partial}{\partial p_{2^{-k}l}}
\Bigr)^2
G_l
=
l^2 \Bigl(p_{2^l}\frac{\partial}{\partial p_l}\Bigr)^2 G_l
+
l^2\Bigl(\frac{\partial}{\partial p_l}\Bigr)^2 G_l
\]
and
\[
l \frac{\partial}{\partial p_l} G_l\Bigm|_{p_l=0}
=
\Bigl(
c_l
+
\sum_{\substack{1\le k\\ 2^k|l}}
2^{-k}l p_{2^{-k}l} \frac{\partial}{\partial p_{2^{-k}l}}
\Bigr)
G_{l-1}.
\]
\end{lem}

\begin{proof}
 From the previous lemma, linearity, and the fact that $X[p_k]$ depends
 only on the variables $p_{2^jk}$, we find that
\[
l C[p_l] \frac{\partial}{\partial p_l} \log(G)
=
c_l
+
\sum_{\substack{1\le k\\2^k|l}}
2^{-k}l p_{2^{-k}l} \frac{\partial}{\partial p_{2^{-k}l}} \log(G),
\]
and thus
\[
l C[p_l] \frac{\partial}{\partial p_l} G
=
\Bigl(
c_l
+
\sum_{\substack{1\le k\\2^k|l}}
2^{-k}l p_{2^{-k}l} \frac{\partial}{\partial p_{2^{-k}l}}
\Bigr) G,
\]
The second differential equation is immediate (since $C[0]=1$); for the
first equation, we have (note that if we set $p_2=p_3=\dots=0$ in $C$, we
obtain the function $\cosh(\arcsinh(p_1))=\sqrt{1+p_1^2}$)
\[
l \sqrt{1+p_l^2} \frac{\partial}{\partial p_l} G_l
=
\Bigl(
c_l
+
\sum_{\substack{1\le k\\2^k|l}}
2^{-k}l p_{2^{-k}l} \frac{\partial}{\partial p_{2^{-k}l}}
\Bigr) G_l.
\]
Since the differential operators on either side commute with each other, we
in fact have
\[
\Bigl(l \sqrt{1+p_l^2} \frac{\partial}{\partial p_l}\Bigr)^2 G_l
=
\bigl(
c_l
+
\sum_{\substack{1\le k\\2^k|l}}
2^{-k}l p_{2^{-k}l} \frac{\partial}{\partial p_{2^{-k}l}}
\Bigr)^2 G_l,
\]
which simplifies to the desired equation.
\end{proof}

\begin{lem}
Suppose $\chi$ is a virtual character of $\Bij$ with cycle index
\[
Z_\chi=\exp(\sum_{l\ge 1} c_l X[p_l])
\]
for some sequence $c_l$ independent of $p_1$, $p_2$,\dots.
Let $\pi$ be a permutation, and for each $m>0$ let $n_m(\pi)$ be the number
of $m$-cycles of $\pi$; also define
\[
o_m(\pi) = \sum_{\substack{1\le k\\ 2^k|m}} 2^{-k}m n_{2^{-k}m}(\pi).
\]
Then
\[
\chi
=
\prod_{l\ge 1} (c_l+o_l)\prod_{0\le i\le n_l-2} (c_l+o_l+l(n_l-2-2i)).
\]
\end{lem}

\begin{proof}
Note that if $n_l=0$, the inner product is over $-1$ terms, and is thus by
standard convention equal to $(c_l+o_l)^{-1}$, so the product over $l$ is
well-defined.

In terms of the cycle index, $\chi$ is given by
\[
\chi(n_1,n_2,\dots)
=
\Bigl(\prod_{1\le l} (l\frac{\partial}{\partial p_l})^{n_l}\Bigr)
Z_\chi\Bigm|_{p_1=p_2=\dots=0},
\]
where we view $\chi$ as a function of the values $n_i(\pi)$.
In particular, the lemma can be interpreted as giving recurrences for the
character; we find that, if $n_{m+1}=n_{m+2}=\dots=0$,
\[
\chi(n_1,\dots,n_m+2,0,0,\dots)
=
((c_m+o_m)^2-(mn_m)^2)
\chi(n_1,\dots,n_m,0,0,\dots)
\]
and similarly from the second differential equation of the lemma,
\[
\chi(n_1,\dots,n_{m-1},1,0,0,\dots)
=
(c_m+o_m)
\chi(n_1,\dots,n_{m-1},0,0,0,\dots).
\]
But then by induction on the sequence $n_i$, in reverse lexicographic
order, the given character formula follows.
\end{proof}

\begin{rem}
The special case $c_1=-c_2=-c_4=-c_8=\cdots=\lambda$, all other $c_i=0$, was
shown in \cite[Thm. 5.7]{CHR}, via a rather different argument.
\end{rem}

In particular, the cycle index of $\Lambda'$ is of this form, and we thus
obtain the following.

\begin{theorem}
Let $\pi$ be a permutation with $n_m(\pi)=n_m$ for $m\ge 1$.  Then
\begin{align}
\chi_{\Lambda'}(\pi;t)
:=&
\Tr(\Lambda'(\pi))\notag\\
=&
\prod_{1\le l}
(\gamma_l(t)+o_l t^{l/2})
\prod_{0\le i\le n_l-2} (\gamma_l(t)+(o_l+l(n_l-2-2i))t^{l/2}),
\notag
\end{align}
where the polynomials $\gamma_l(t)$ are given by the expression
\[
\gamma_l(t) = \sum_{\substack{k|l\\\text{$k$ odd}}} \mu(k) t^{l(1-1/k)/2},
\]
where $\mu$ is the M\"obius function.
\end{theorem}

\begin{proof}
This is equivalent to the claim
\[
Z_{\Lambda'} = \exp(\sum_{1\le l} t^{-l/2}\gamma_l(t) X[t^{l/2} p_l]/l)
\]
since then we can apply the lemma to $Z_{\Lambda'}[t^{-1/2} p_1]$.

The claimed expression for $Z_{\Lambda'}$ is easily obtained by expanding
\begin{align}
\sum_{1\le l} t^{-l/2} \Arcsinh^\sim[t^{l/2} p_l]/l
&=
\sum_{1\le l} t^{-l/2} \sum_{\text{$m$ odd}} \mu(m) X[t^{lm/2} p_{lm}]/lm
\notag\\
&=
\sum_{1\le l} 
\sum_{\substack{m|l\\\text{$m$ odd}}}
t^{-(l/m)/2} \mu(m) X[t^{l/2} p_l]/l.
\notag
\end{align}
\end{proof}

\begin{cor}\label{cor:char}
Theorem \ref{thm:main} holds whenever $n_1\ge 0$.
\end{cor}

\begin{proof}
If $n_1\ge 0$, or in other words if $\pi$ has a fixed point (so WLOG
$\pi(n)=n$), then this follows immediately from the isomorphism between
$\Lambda'(\{1,2,\dots,n-1\})$ and $H^*(M_n,\QQ)$.
\end{proof}

\section{Functoriality}

In fact, as we will see, the character formula continues to hold even if
$\pi$ has no fixed point (so $n_1=-1$).  The key idea is that although we
have so far only considered $\Lambda$ as a functor on $\Bij$ (or more
precisely on the category of {\em nonempty} finite sets and bijections), it
actually extends to a functor on the full category $\Fin^+$ of
nonempty sets.

For a nonempty finite set $S$, let $\Lambda(S)$ denote the algebra
isomorphic to $\Lambda_n$ with generators $\omega_{ijkl}$ for $i,j,k,l\in
S$.  This extends to a functor $\Lambda:\Fin^+\to \GrAlg$ as follows.  If
$f:S\to T$ is an arbitrary function, we define
\[
\Lambda(f)(\omega_{ijkl})
=
\omega_{f(i)f(j)f(k)f(l)},
\]
where $\omega_{ijkl}:=0$ if any two indices are equal.  Since this
convention makes the defining relations of $\Lambda$ hold even if some
indices coincide, we indeed obtain a homomorphism.

This has important consequences for the $S_n$-module structure, as the
irreducible representions of the category $\Fin^+$ are easily determined
(and defined over $\QQ$).  The irreducible representation theory of
$\Fin^+$ is determined by the irreducible representation theory of the
``transformation semigroup'' (the semigroup of functions from a finite set
to itself).  Thus from results of \cite{Put}, we immediately have the
following (compare chapter 8 of \cite{Ner}).

\begin{theorem}
Let $R$ be an irreducible complex representation of $\Fin^+$.  Then
precisely one of the following two statements holds for $R$.
\begin{itemize}
\item[1] There exists a nonnegative integer $k$ such that $R(\{1,2,\dots,
  n\})$ is $\binom{n-1}{k}$-dimensional, with $S_n$-character with label
  $(n-k) 1^k$ for $n>k$.
\item[2] There exists a partition $\lambda$ not of the form $1^k$ such that
  each $S_n$-module $R(\{1,2,\dots,n\})$ is induced from the
  $S_{|\lambda|}\times S_{n-|\lambda|}$-module in which $S_{n-|\lambda|}$ acts
  trivially and $S_{|\lambda|}$ acts as the representation $\lambda$.
\end{itemize}
In particular, we can choose a basis of each $R(S)$ such that all matrix
coefficients are rational.
\end{theorem}

\begin{remark}
Note, however, that the transformation semigroup does not have
finite representation type, and thus the full representation theory of
$\Fin^+$ is wild.
\end{remark}

If $R$ is an irreducible representation of $\Fin^+$ with cycle
index $Z_R$, then in the first case we have
\[
Z_R = (-1)^k+(e_k-e_{k-1}+e_{k-2}+\dots)\Exp,
\]
where $e_k$ is the cycle index of the sign representation of $S_k$, while
in the second case we have $Z_R = s_\lambda \Exp$, where $s_\lambda$ is a
Schur function (the cycle index of the irreducible representation indexed
by $\lambda$).

\begin{cor}
If $R$ is a representation of $\Fin^+$ with cycle index $Z_R$, such that
$\dim(R(S))=O(|S|^l)$ for some integer $l\ge 0$ (``polynomial growth''),
then there exists a unique constant $C_R$ such that $\Exp^{-1}(C_R+Z_R)$ is
a symmetric function of $\deg_p$ degree at most $l$.
\end{cor}

\begin{proof}
Since $\dim(R(S))=O(|S|^l)$, the same must be true for the irreducible
constituents of $R$, which must therefore satisfy $k\le l$ or
$|\lambda|\le l$, as appropriate.  The result follows.
\end{proof}

If $R$ is such a representation (or more generally, a graded representation
in which each homogeneous component has polynomial growth), we will call
$C_R+Z_R$ the {\em extended} cycle index of $R$, and denote it by $Z^+_R$.
Note in particular that $\dim(R(S))$ is polynomial in $|S|$, with constant
term $C_R$.

\begin{cor}
The coefficient of $t^k$ in $\Exp^{-1}Z^+_\Lambda$ is a symmetric function
of degree at most $3k$.
\end{cor}

\begin{proof}
Indeed, the formula for the Poincar\'e series of $\Lambda$ implies that
the degree $k$ component of $\Lambda(S)$ has dimension $O(n^{3k})$.
\end{proof}

Now, it follows easily from the fact that the $S_{n-1}$-module
$\Lambda'_{n-1}$ is the restriction of the $S_n$-module $\Lambda_n$ that
\[
Z_{\Lambda'} = \frac{\partial}{\partial p_1} Z_\Lambda.
\]
But this together with the corollary is enough to uniquely determine
$Z^+_\Lambda$.  Indeed, in general, if
\[
C_R+Z_R = f\Exp
\]
for some symmetric function $f$ of finite degree, then
\[
\frac{\partial}{\partial p_1} Z_R = (f+\frac{\partial}{\partial p_1}
f)\Exp.
\]
If we write $f$ as a polynomial in $p_1$, we can then solve for its
coefficients in order starting with the highest degree term; in other
words, the (extended) cycle index of any $\Fin^+$ representation with
polynomial growth is uniquely determined by the cycle index of its
restriction to point stabilizers.

In our case, we can explicitly solve the corresponding differential
equation.

\begin{thm}\label{thm:ext_cycle}
The extended cycle index of the $\Fin^+$ representation $\Lambda$
is given by
\[
Z^+_{\Lambda}
=
\frac{-p_1 t + \Cosh^\sim[\Arcsinh^\sim[t^{1/2}p_1]]}{1-t}
\Exp[t^{-1/2} \Arcsinh^\sim[t^{1/2} p_1]].
\]
In particular, $C_\Lambda = 1/(1-t)$.
\end{thm}

\begin{proof}
To prove the theorem, we need simply verify that the above
expression differentiates to $Z_{\Lambda'}$ and that if we divide by $\Exp$ the
coefficient of $t^k$ is of bounded degree.

If we divide the above expresion by $\Exp$, we obtain
\[
\frac{-p_1 t + \Cosh^\sim[\Arcsinh^\sim[t^{1/2}p_1]]}{1-t}
\Exp[t^{-1/2} \Arcsinh^\sim[t^{1/2} p_1]-p_1].
\]
Now, if $f$ and $g$ are symmetric functions satisfying the bounded degree
condition, then so are $f+g$ and $fg$; if moreover $g$ has constant term 0
as a series in $t$, then $f[g]$ has bounded degree coefficients.  The
second condition follows.

 From the identities
\begin{align}
\frac{\partial}{\partial p_1} \Exp[f]
&= (\frac{\partial}{\partial p_1} f)\Exp[f],\notag\\
\frac{\partial}{\partial p_1} \Sinh^\sim[f]
&= (\frac{\partial}{\partial p_1} f)\Cosh^\sim[f],\notag\\
\frac{\partial}{\partial p_1} \Cosh^\sim[f]
&= (\frac{\partial}{\partial p_1} f)\Sinh^\sim[f],\notag
\end{align}
we find, differentiating the defining equation for $\Arcsinh^\sim$, that
\[
(\frac{\partial}{\partial p_1} \Arcsinh^\sim) \Cosh^\sim[\Arcsinh^\sim] = 1
\]
and can then immediately verify that $Z^+_\Lambda$ differentiates as required.
\end{proof}

\begin{rems}
1. The above formula was guessed via the corresponding formula for the
(super) Poincar\'e series (i.e., setting $p_1=u$, $p_2=p_3=\dots=0$):
\[
\frac{-ut+\cosh(\arcsinh(u\sqrt{t}))}{1-t}
\exp(\arcsinh(u\sqrt{t})/\sqrt{t}).
\]

2.  Getzler \cite{Getz} gave the complex analogue, as follows.  If we
subtract $p_1$ from the cycle index of $H^*(\overline{M_{0,|S|}},\CC)$ we
get the unique solution $Z$ of
\begin{align}
Z &= p_1 \frac{\partial}{\partial p_1} Z - F[\frac{\partial}{\partial p_1}
  Z],\notag\\
F &= p_1 \frac{\partial}{\partial p_1} F - Z[\frac{\partial}{\partial p_1}
  F]\notag
\end{align}
where
\[
F = \frac{\Exp[(1+t^2)\Log[1+p_1]]-(1+p_1)(1+t^2 p_1)}{t^2-t^6}
+\frac{h_2}{1+t^2}
+\frac{p_1^2-p_2}{2}
\]
(essentially the cycle index for $M_{0,|S|}(\CC)$; note that the formula
given in \cite{Getz} is slightly incorrect, but the correct formula follows
from the results of \cite{KL}).  More precisely, the
full cycle index can be expressed as
\[
p_1 \Exp[C] - F[\Exp[C]-1],
\]
where $C$ is as in Theorem \ref{thm:cyc_ind_C} above.  The equations
relating $Z$ and $F$ constitute an involutory transformation (the
``Legendre transform'') integrating the plethystic inverse.  With this in
mind, we note that the Legendre transform of $Z^+_\lambda-(1-t)^{-1}-p_1$ is
\[
\frac{(1+p_1)(\Cosh^\sim[t^{1/2}\Log[1+p_1]]-t^{-1/2}\Sinh^\sim[t^{1/2}\Log[1+p_1]])-1}{1-t},
\]
which again presumably has a cohomological interpretation.
\end{rems}

This also allows us to prove the remaining cases of Theorem \ref{thm:main}.

\begin{proof}
The point is that if $R$ is any representation of $\Fin^+$ with polynomial
growth, then the character of $R$ depends polynomially on the numbers $n_i$
of $i$-cycles (since this holds for irreducibles).  In particular, for
$\Lambda$ we may thus extrapolate to the case with no fixed points.
\end{proof}

If we set $t=1$ in the formula for the graded character, we obtain the
Euler character of $M_n$.  This is straightforward
except in the case $n_1=-1$, when we have a factor $1/(1-t)$ that must be
cancelled.  We obtain the following result.

\begin{cor}\label{cor:euler}
The Euler character $\chi_E$ of $M_n$ at the permutation $\pi\in S_n$ is
nonzero if and only if one of the following (disjoint) conditions is
satisfied.  We suppose $\pi$ has $n_1+1$ fixed points and $n_l$ $l$-cycles
for $l\ge 2$.
\begin{itemize}
\item[1.] $\pi$ has a fixed point.  Then $\pi$ has order a power of 2 and
  there exists $k\ge 0$ with $n_1=n_2=\dots n_{2^{k-1}}=1$, $n_{2^k}$ even.
  In this case,
\[
\chi_E(\pi)
=
2^{k(k-1)/2}
\prod_{k\le j} (1+o_{2^j})
\prod_{0\le i\le n_{2^j}-2} (1+o_{2^j}+2^j(n_{2^j}-2-2i)).
\]
\item[2.] $\pi$ has no fixed points.  Then there exists a nonnegative integer
  $d>1$ such that $n_d$ is odd, and every cycle of $\pi$ has length $2^j d$
  for some $j$.  In this case,
\begin{align}
\chi_E(\pi)
={}&
\bigl(\sum_{\text{odd }k|d} \frac{\mu(k)d}{2k}\bigr)
\prod_{0\le i\le n_d-2}d(n_d-2-2i)\notag\\
&\prod_{1\le j} (1+o_{2^jd})
\prod_{0\le i\le n_{2^jd}-2} (1+o_{2^jd}+2^jd(n_{2^jd}-2-2i)).
\notag
\end{align}
\end{itemize}
\end{cor}

\begin{proof}
If $\pi$ has a fixed point, then we may simply set $t=1$ in the formula for
the graded character:
\[
\prod_{1\le l}
(\gamma_l(1)+o_l)
\prod_{0\le i\le n_l-2} (\gamma_l(1)+(o_l+l(n_l-2-2i))),
\]
where we recall that $\gamma_l(1)=0$ unless $l$ is a power of 2, in which
case $\gamma_l(1)=1$.  Suppose $\pi$ did not have order a power of 2; then
in particular it would have a cycle of length not a power of 2.  Let $d$ be
the length of the shortest such cycle.  Then $o_d=\gamma_d(1)=0$, and the
contribution of the $l=d$ factor of the above product is 0.

Similarly, let $k$ be the smallest integer such that $n_{2^k}\ne 1$.  Then
$o_{2^k}=2^k-1$, and the contribution for $l=2^k$ is
\[
2^k
\prod_{0\le i\le n_{2^k}-2} (2^k+2^k(n_{2^k}-2-2i)).
\]
If $n_{2^k}$ were odd, then the factor for $i=(n_{2^k}-1)/2$ would make the
product 0, and thus $n_{2^k}$ must be even.  (In particular, it follows
that $n/2^k$ is odd.)  The above formula for the Euler character is
then straightforward.  Since $n_{2^k}$ is even, it follows that $o_{2^j}/2^k$ is
odd for all $j$, and thus none of the remaining factors can vanish.

Now, suppose $\pi$ has no fixed points.  In this case, the contribution for
$l=1$ to the graded character is a factor $1/(1-t)$, and thus rather than
avoid all factors that vanish for $t=1$, we must have exactly one vanishing
factor.  Suppose $d$ is the length of the shortest cycle of $\pi$.  Then
$\gamma_d(1)+o_d=0$ (since either $d$ is a power of 2, with $o_d=n_1=-1$,
or $d$ is not a power of 2, and $o_d=0$), and thus it provides that
vanishing factor (and provides more than one unless $n_d$ is odd; in
particular $n/d$ must be odd).  If
there were a cycle of any length not of the form $2^j d$, the shortest such
cycle would provide {\em another} vanishing factor.  The above formulae for
the Euler character are again straightforward.
\end{proof}

\begin{rems}
1. For the Euler characteristic itself, the above criterion translates to
the statement that $\chi=0$ unless $n$ is odd (i.e., $n_1=n-1$ is even), in
which case the Euler characteristic is
\[
\prod_{0\le i\le n-3} (1+(n-3-2i))
=
\prod_{0\le i\le \lfloor (n-3)/2\rfloor} (1-(n-3-2i)^2),
\]
agreeing with the calculations of \cite[Thm. 3.2.3]{Dev} and \cite{Gai}.

2. It should be possible to prove this directly by studying the fixed point
set of the action of $\pi$ on $M_n$.
\end{rems}


\begin{thebibliography}{9}

\bibitem{Bj} A. Bj\"orner,
Shellable and {C}ohen-{M}acaulay partially order sets.
Trans. Amer. Math. Soc. 260 (1980), no. 1, 159--183.

\bibitem{CHR} A. R. Calderbank, P. Hanlon, R. W. Robinson,
Partitions into even and odd block size and some unusual characters of the
symmetric groups.
Proc. London Math. Soc. (3) 53 (1986), no. 2, 288--320.

\bibitem{Dev} S. Devadoss,
Tessellations of Moduli Spaces and the Mosaic Operad,
Contemporary Mathematics 239 (1999), p. 91-114

\bibitem{EHKR} P. Etingof, A. Henriques, J. Kamnitzer, E. Rains.
The cohomology ring of the real locus of the moduli space of stable curves
of genus 0 with marked points.
arXiv:math.AT/0507514.

\bibitem{Gai}
G. Gaiffi, Real structures of models of arrangements.  Int. Math. Res. Not.  200
4,  no. 64, 3439--3467.

\bibitem{Getz}
E. Getzler, The semi-classical approximation for modular operads.
Commun. Math. Phys. 194 (1998), 481--492.

\bibitem{GK}
V. Ginzburg, M. Kapranov.
Koszul duality for operads.
Duke Math. J.  76  (1994),  no. 1, 203--272.

\bibitem{Hen}
A. Henderson,
Representations of wreath products on cohomology of {D}e
{C}oncini-{P}rocesi compactifications.
Int. Math. Res. Not. 2004, no. 20, 983--1021.

\bibitem{KL}
M. Kisin, G. I. Lehrer,
Equivariant Poincar\'e polynomials and counting points over finite fields.
J. Algebra 247 (2002), no. 2, 435--451.

\bibitem{Leh}
G. I. Lehrer,
On the Poincar\'e series associated with Coxeter group actions on
complements of hyperplanes,
J. London Math. Soc. (2) 36 (1987), 275--294.

\bibitem{Ner} Yu. A. Neretin,
Categories of symmetries and infinite-dimensional groups,
London Mathematical Society Monographs (NS), 16. Oxford University Press,
New York, 1996.

\bibitem{Put} M. S. Putcha,
Complex representations of finite monoids, Proc. London Math. Soc. (3), 73
(1997), no. 3, 623--641.

\bibitem{Yuz} S. Yuzvinsky,
Cohomology bases for the De Concini-Procesi models of hyperplane
arrangements and sums over trees.
Invent. Math. 127 (1997), no. 2, 319--335.

\end{thebibliography}
\end{document}